\documentclass[a4paper,11pt,reqno]{amsart}
\usepackage{mathrsfs}
 \usepackage[all]{xy}
\usepackage{amsmath,amssymb,amscd,bbm,amsthm,mathrsfs}
\usepackage{color,xcolor}
\usepackage{graphicx}
\usepackage{manfnt}
\newtheorem{thm}{Theorem}[section]
\newtheorem{lem}{Lemma}[section]
\newtheorem{cor}{Corollary}[section]
\newtheorem{con}{Conjecture}[section]

\newtheorem{rem}{Remark}[section]

\begin{document}
\numberwithin{equation}{section}
 \title[Quasi-positive mixed curvature]{Quasi-positive mixed curvature, vanishing theorems, and rational connectedness}
\author{Kai Tang}
\address{Kai Tang. School of Mathematical Sciences, Zhejiang Normal University, Jinhua, Zhejiang, 321004, China} \email{{kaitang001@zjnu.edu.cn}}
\thanks{\text{Foundation item:} Supported by National Natural Science
Foundation of China (Grant No.12001490).}
\begin{abstract}
In this paper, we consider {\em mixed curvature} $\mathcal{C}_{a,b}$, which is  a convex combination of Ricci curvature and holomorphic sectional curvature introduced by Chu-Lee-Tam \cite{CLT}. We prove that if a compact complex manifold $M$ admits a K\"{a}hler metric with quasi-positive mixed curvature and $3a+2b\geq0$, then it is projective. If $a,b\geq0$, then $M$ is rationally connected. As a corollary, the same result holds for $k$-Ricci curvature. We also show that any compact K\"{a}hler manifold with quasi-positive 2-scalar curvature is projective. Lastly, we  generalize the result to the Hermitian case. In particular, any compact Hermitian threefold with quasi-positive real bisectional curvature have vanishing Hodge number $h^{2,0}$. Furthermore, if it is K\"{a}hlerian, then it is projective.
\end{abstract}

 \maketitle

\section{Introduction}
\subsection{Background}
For compact K\"ahler manifolds with positive holomorphic sectional curvature or positive Ricci curvature, Yau proposed the following conjecture.
\begin{con}[\cite{Yau}, Problem 47]
	Let $(M,g)$ be a compact K\"ahler manifold with positive holomorphic sectional curvature or positive Ricci curvature, then $M$ is projective and rationally connected.
\end{con}
A projective manifold $M$ is called rationally connected if any two points of $M$ can be joined by a rational curve. Thanks to the celebrated Calabi-Yau theorem \cite{Yau78}, we know that a compact K\"ahler manifold $M$ has a Hermitian metric with positive first Chern-Ricci curvature
if and only if $M$ is Fano. The Campana \cite{Cam} and Kollar-Miyaoka-Mori \cite{KMM} showed that a K\"ahler manifold with positive Ricci curvature is rationally connected. Then this result was generalized to the quasi-positive first or second Chern-Ricci curvature by X. Yang in \cite{Yang2020}. In \cite{LZZ}, Li-Zhang-Zhang established the equivalence between rational connectedness and the quasi-positivity of mean curvature for compact K\"ahler manifolds.

Recently, X. Yang \cite{Yang2018} proved that compact K\"ahler manifolds with positive holomorphic sectional curvature are all projective. Hence by the previous work of Heier-Wong \cite{HW} that any projective manifold with quasi-positive holomorphic sectional curvature must be rationally connected, Yau's conjecture on positive holomorphic sectional curvature was solved. In \cite{Mat1,Mat2}, Matsumura established the structure theorems for a projective manifold with nonnegative holomorphic sectional curvature. In particular, he strengthened Heier-Wong's result \cite{Mat2}: if a projective manifold $M$ admits a K\"{a}hler metric $g$ with nonnegative holomorphic sectional curvature and $(M, g)$ has no nonzero truly flat tangent vector at some point (which is satisfied when the holomorphic sectional curvature is quasi-positive), $M$ is rationally connected. Motivated by these results, Yau's conjecture on positive holomorphic sectional curvature was conjectured \cite{Mat3} by weakening the assumption to quasi-positivity, which was recently solved by S. Zhang and X. Zhang in \cite{ZZ}. The proof crucially relied on a new Bochner-type integral equality.

In this paper, we shall employ the methods in \cite{ZZ} (Zhang-Zhang) and \cite{HW} (Heier-Wong) to consider the quasi-positive case for some other interesting curvature conditions.

\subsection{Mixed curvature $\mathcal{C}_{a,b}$}
For a K\"{a}hler manifold $(M^{n},g)$, the first Ricci curvature $Ric$  of the Chern connection is defined by
\begin{align}
Ric=\sqrt{-1}R_{i\overline{j}}dz^{i}\wedge d\overline{z}^{j}=\sqrt{-1}(g^{k\overline{l}}R_{i\overline{j}k\overline{l}})dz^{i}\wedge d\overline{z}^{j}=-\sqrt{-1}\partial\overline{\partial}\log \det g\,. \nonumber
\end{align}
where it is a $(1,1)$-form representing the first Chern class $c_{1}(M)$. The holomorphic sectional curvature $H$ is defined by
\begin{align}
H(X)=R(X,\overline{X},X,\overline{X})/|X|^{4} \nonumber \,\,,
\end{align}
for a (1, 0)-tangent vector $X\in T^{1,0}M$. Given the perplexing relationship between the holomorphic sectional curvature and
the Ricci curvature on K\"{a}hler manifolds,  one can attempt to interpolate between these curvatures
by considering the following curvature constraint: In \cite{CLT} (also see \cite{BT}), Chu-Lee-Tam
introduced, for $a,b\in\mathbb{R}$, the {\em mixed curvature}
 \begin{align}
\mathcal{C}_{a,b}(X)=\frac{a}{|X|^{2}_{g}}Ric(X,\bar{X})+b H(X)\nonumber
 \end{align}
for a (1, 0)-tangent vector $X\in T^{1,0}M$. We shall use the mixed curvature to proceed the following curvature conditions: $\mathcal{C}_{1,0}$ is the Ricci curvature, $\mathcal{C}_{0,1}$ is the holomorphic sectional curvature, and $\mathcal{C}_{1,-1}$  is the orthogonal Ricci curvature $Ric^{\perp}(X,\overline{X})$, a curvature condition defined  as the difference between Ricci curvature and holomorphic sectional curvature by Ni-Zheng \cite{NZ1}. In addition, $\mathcal{C}_{k-1,n-k}$ is closely related to the $k$-$Ricci$ curvature. In \cite{CLT, BT},
it was shown that any compact K\"{a}hler manifold with positive mixed curvature, which constants satisfy $a>0 $ and $3a+2b\geq0$, must be projective and simply connected. The first main result of our paper is the following:
\begin{thm}\label{thm1.1} Let $(M^{n},g)$ be a compact K\"{a}hler manifold with quasi-positive mixed curvature $\mathcal{C}_{a,b}$. Then the following statements holds:
\begin{itemize}
\item[(1)] If  $a\geq0$ and $3a+2b\geq0$, then $h^{2,0}=0$. In particularly, $M$ is projective.
\item[(2)] If $a,b\geq0$, then $M$ is rationally connected.
\end{itemize}
\end{thm}
As a direct application, we obtain that any compact K\"{a}hler manifold with quasi-positive orthogonal Ricci curvature must be projective. In an attempt to generalize the hyperbolicity of Kobayashi to the $k$-hyperbolicity, Ni \cite{N1} introduced the concept of {\em $k$-Ricci curvature}. Given a compact K\"{a}hler manifold $(M^{n}, g)$. The $k$-Ricci curvature $Ric_{k}$ ($1\leq k\leq n$) is defined as the Chern Ricci curvature of the $k$-dimensional holomorphic subspaces of the holomorphic tangent bundle $T^{1,0}M$. It was proved by Ni \cite{N2} that a compact K\"{a}hler manifold with positive $k$-Ricci curvature  must be projective and rationally connected. As a corollary of Theorem \ref{thm1.1}, we also generalize the Ni's result to the quasi-positive case (see Remark \ref{rem1.3} for the proof).
\begin{cor}\label{1.3} Let $(M^{n}, g)$ be a  compact K\"{a}hler manifold  with quasi-positive $k$-Ricci curvature. Then $M$ must be projective and rationally connected.
\end{cor}

\subsection{$2$-scalar curvature $S_2$}
Another interesting curvature condition is {\em $k$-scalar curvature}, which is studied in \cite{NZ2}. Let $(M^{n}, g)$ be a compact Hermitian manifold. The $k$-scalar curvature $S_{k}$ ($1\leq k\leq n$) is defined as the scalar curvature of the $k$-dimensional holomorphic subspaces of the holomorphic tangent bundle $T^{1,0}M$. If $g$ is K\"{a}hler, by Berger's lemma \cite{Berger} one can define the $k$-scalar curvature of $g$ as
\begin{align}
S_{k}(p,\Sigma)=\frac{k(k+1)}{2 Vol(\mathbb{S}^{2k-1})}\int_{|Z|=1,Z\in\Sigma}H(Z)d\theta(Z) \nonumber
\end{align}
for any $p\in M$ and any $k$-dimension subspace $\Sigma\in T^{1,0}M$. When $g$ is K\"{a}hler metric, the positivity of the holomorphic sectional curvature implies the positivity of
the 2-scalar curvature, and the positivity of $S_{k}$ implies the positivity of $S_{l}$ if $k\leq l$.

Ni-Zheng \cite{NZ2} showed that any compact K\"{a}hler manifold with positive $2$-scalar curvature must be projective. However, in generally, the projectivity of $M$ cannot be implied by the positivity of $S_k$ for $k\geq 3$ (taking the product of a very positive $\mathbb{P}^1$ and a non-projective torus of complex dimension $2$). The dimension of the fiber of the MRC fibration of $M$ is called the {\em rational dimension} of $M$, and is denoted
by $rd(M)$. Heier-Wong \cite{HW} proved that a projective manifold with quasi-positive $S_k$ satisfies $rd(M)\geq n-(k-1)$. The second main result of this paper is the following theorem, which generalizes Ni-Zheng's result to the quasi-positive case, and also extends Heier-Wong's result to compact K\"ahler manifolds for $k=2$ (see Remark \ref{rem1.3} for the proof).
\begin{thm}\label{1.4} If $(M^{n}, g)$ is a compact K\"{a}hler manifold  with quasi-positive $2$-scalar curvature, then $h^{2,0}=0$. In particular, $M$ must be projective and $rd(M)\geq n-1$.
\end{thm}

\subsection{Some special case for the Hermitian category}
For the Hermitian case, if metric is only assumed to be Hermitian metric with  positive or quasi-positive  holomorphic sectional curvature and the K\"{a}hlerity is a priori unknown, then it is difficult to ensure the projectivity of compact Hermitian manifolds even under a slightly stronger curvature condition, {\em real bisectional curvature} introduced by Yang-Zheng \cite{YZ}. They proposed the following generalization of Yau's conjecture for Hermitian case.
\begin{con} Let $(M,g)$ be a Hermitian manifold with positive (quasi-positive) real bisectional curvature. Then $M$ must be a projective and rationally connected.
\end{con}
Under the entra assumption that $M$ is K\"{a}hlerian, X. Yang \cite{Yang20201} also conjectured that $M$ is projective and rationally connected if $M$ adimts a Hermitian metric with quasi-positive holomorphic sectional curvature. He showed that it is  true if surface $M$ adimts a Hermitian metric with positive holomorphic sectional curvature. In \cite{Tang1, Tang2}, the author proved the the projectivity of  compact Hermitian manifold $M$ with positive real bisectional curvature under the additional assumption that $M$ is K\"{a}hlerian. Furthermore, such manifold is rationally connected, which was proved by Zhang-Zhang \cite{ZZ}. For quasi-positive case, it is difficult to prove that the Hodge number $h^{2,0}$ equals 0. But we can also obtain some meaningful results.
\begin{thm}\label{1.6} Let $(M^{3}, g)$ be a compact Hermitian threefold with non-negative real bisectional curvature. Then the following statements holds:
\begin{itemize}
\item[(1)] For any holomorphic $(2,0)$-form $\eta$, we have that $|\eta|_{g}^{2}\equiv C$ for some constant $C\geq0$.
\item[(2)] If the real bisectional curvature is quasi-positive, then $H^{2,0}_{\bar{\partial}}(M)=0$. In particular, if
$M^{3}$ is K\"{a}hlerian, then it is projective.
\end{itemize}
\end{thm}
If the dimension is greater than 3,  we have the following result.
\begin{thm}\label{1.7} Let $(M^{n}, g)$ be a compact Hermitian manifold with non-negative real bisectional curvature. If $M$ is K\"{a}hlerian and $\partial\bar{\partial}\omega^{n-3}=0$, then the following statements holds:
\begin{itemize}
\item[(1)] For any holomorphic $(2,0)$-form $\eta$, we have that $|\eta|_{g}^{2}\equiv C$ for some constant $C\geq0$.
\item[(2)] If the real bisectional curvature is quasi-positive, then $H^{2,0}_{\bar{\partial}}(M)=0$. In particular, $M$ is projective.
\end{itemize}
\end{thm}
A special case of Theorem \ref{1.7} is
\begin{cor}\label{1.8} Let $M^{4}$ be a compact K\"{a}hler  four-manifold with a pluriclosed metric of quasi-positive real bisectional curvature. Then $M$ is projective.
\end{cor}
\begin{rem} If we replace real bisectional curvature with 2-scalar curvature, the results of Theorem \ref{1.6}  still hold. For the Hermitian case, the positivity of real bisectional curvature can not imply the positivity of 2-scalar curvature.
\end{rem}

To conclude the introduction, we emphasize the key ingredient of this work. The proofs combine two distinct approaches: Zhang-Zhang's method \cite{ZZ} to establish projectivity and Heier-Wong's technique \cite{HW} to prove rational connectedness. Thanks to Berger's averaging trick, we are able to carry out some new computations that utilize linear combinations of the Bochner integral formulas for both holomorphic sectional curvature and Ricci curvature.
This method should be able to deal with other interesting curvature conditions.

This paper is organized as follows. In section 2, we provide some basic knowledge which will be used in our proofs. In section 3, we prove  Theorem 1.1, Corollary 1.1, Theorem 1.2. In section 4, we prove  Theorem 1.3 and Theorem 1.4.

\vspace{0.5cm}
\noindent\textbf{Acknowledgement.} The author is grateful to Professor Fangyang Zheng for constant encouragement and support.  We are grateful to Shiyu Zhang  for useful discussions on the results related to mixed curvature and for suggestions which improved the readability of the paper. We also thank the referees for their helpful comments which enhance the presentation.

\section{Preminaries}
Let $(M,g)$ be a Hermitian manifold of dimension $n$. Let $\alpha$ be a real $(1,1)$-form and $\eta$ be a $(p,0)$-form on $M$. As in \cite{ZZ}, we can define a real semi-positive $(1,1)$-form $\beta$ associated to $\eta$ by
 \begin{align}
\beta=\Lambda^{p-1}((\sqrt{-1})^{p^{2}}\frac{\eta\wedge \bar{\eta}}{p!}) \nonumber
 \end{align}
 where $\Lambda$ is the dual Lefschetz operator with respect to $\omega_{g}$. In a local orthonormal coordinates at $x_{0}$, we write $\eta=\eta_{I_p}dz^{I_p}$, $\omega=\sqrt{-1}dz^{i}\wedge d\bar{z}^{i}$ and $\beta=\sqrt{-1}\beta_{i\bar{i}}dz^i\wedge d\bar{z}^i$, then
 $$\beta_{i\bar{i}}=p!\sum\limits_{I_{p-1}}\eta_{iI_{p-1}}\overline{\eta_{iI_{p-1}}}$$
 Hence $tr_{\omega}\beta=|\eta|_{g}^{2}$. We have the following formula.

\begin{lem}[\cite{ZZ}, lemma 3.2]
	\begin{align}\label{2.1}
(\sqrt{-1})^{p^{2}}\alpha \wedge \eta \wedge \bar{\eta} \wedge \frac{\omega^{n-p-1}}{(n-p-1)!}=[tr_{\omega}\alpha\cdot|\eta|_{g}^{2}-
p\langle\alpha,\beta\rangle_{g}]\frac{\omega^{n}}{n!}
 \end{align}
\end{lem}
\begin{lem}[\cite{ZZ}, formula (3.14)] Let $\eta$ be a holomorphic $(p,0)$-form. Then
  \begin{align}\label{2.2}
  	\sqrt{-1}\partial_u\bar{\partial}_v|\eta|^{2}_{g}=\langle D_u'\eta,\overline{D_v'\eta} \rangle+p\sum\limits_{i=1}^nR_{u\bar{v}i\bar{i}}\beta_{i\bar{i}}
  \end{align}
\end{lem}

Another useful result is the following
\begin{lem}[\cite{CLT}, lemma 2.2]
Let $(M^{n},g)$ be a K\"{a}hler manifold with $Ric_{k}(X,\bar{X})\geq (k+1)\sigma|X|^{2}$. Then the following estimates holds:
  \begin{align}\label{2.3}
(k-1)|X|^{2}Ric(X,\bar{X})+(n-k)R(X,\bar{X},X,\bar{X})\geq (n-1)(k+1)\sigma|X|^{4}
  \end{align}	
\end{lem}

We further state Berger's lemma (see \cite{Berger} or \cite{NZ2}), which plays a key role in the computation of mixed curvature.
\begin{lem}[Berger]
Let $(M^{n},g)$ be a K\"{a}hler manifold. If $S(p)=\sum_{i,j=1}^{n}R(e_{i},\overline{e}_{i},e_{j},\overline{e}_{j})$, where $\{e_{i}\}$ is a unitary
basis of $T_{p}^{1,0}M$, denotes the scalar curvature of $M$, then
  \begin{align}\label{2.4}
S(p)=\frac{n(n+1)}{2 Vol(\mathbb{S}^{2n-1})}\int_{|Z|=1,Z\in T_{p}^{1,0}M}H(Z)d\theta(Z)
  \end{align}	
\end{lem}
It is also easy to check that
\begin{align}\label{2.5}
S(p)=\frac{n}{Vol(\mathbb{S}^{2n-1})}\int_{|Z|=1,Z\in T_{p}^{1,0}M}Ric(Z,\overline{Z})d\theta(Z)
  \end{align}

 Finally, let us  recall the  concept of {\em real bisectional curvature}.  Let $(X^{n},g)$ be a Hermitian manifold. Denote by $R$ the curvature tensor of the Chern connection. For $p\in X$, let $e=\{e_{1},\cdot\cdot\cdot,e_{n}\}$ be a unitary tangent frame at $p$, and let $a=\{a_{1},\cdot\cdot\cdot,a_{n}\}$ be non-negative constants with $|a|^{2}=a_{1}^{2}+\cdot\cdot\cdot+a_{n}^{2}>0$. Define the $real$ $bisectional$ $curvature$ of g by
\begin{align}
B_{g}(e,a)=\frac{1}{|a|^{2}}\sum_{i,j=1}^{n}R_{i\overline{i}j\overline{j}}a_{i}a_{j}.
\end{align}
When the metric is K\"ahler, this curvature is the same as the holomorphic sectional curvature $H$, while when the metric is not K\"ahler, the curvature condition is slightly stronger than $H$ at least algebraically. For a more detailed discussion of this, we refer the readers to \cite{YZ}.
	
\section{Mixed curvature}
In this section, we prove Theorem \ref{thm1.1}, and the proof of the Corollary is essentially the same.
\subsection{Projectivity}

To prove the projectivity, we only need to show that $h^{2,0}=0$. Assume that $M$ admits a nonzero holomorphic $(2,0)$-form
$\sigma=\sigma_{ij}dz^{i}\wedge dz^{j}$. Let $k$ be the largest integer such that $\sigma^{k}\neq0$ but $\sigma^{k+1}=0$.

Now, let $x\in M$ be a fixed point such that $\sigma^{k}\neq0$, we can choose local holomorphic orthonormal coordinates such that
 \begin{align*}
\sigma=\lambda_{1}e^{1}\wedge e^{2}+\cdot\cdot\cdot+\lambda_{k}e^{2k-1}\wedge e^{2k}
 \end{align*}
 Let $\eta=\sigma^{k}$, so it is a holomorphic $(2k,0)$-form. Then
 \begin{align*}
\eta=k!\lambda_{1}\cdot\cdot\cdot\lambda_{k}e^{1}\wedge\cdot\cdot\cdot\wedge e^{2k}
 \end{align*}
 Let $\beta$ be the $(1,1)$-form associated to $\eta$, we have that
 \begin{align*}
\beta&=\frac{(k!\lambda_{1}\cdot\cdot\cdot\lambda_{k})^{2}}{2k}(e^{1}\wedge\bar{e}^{1}+\cdot\cdot\cdot e^{2k}\wedge\bar{e}^{2k})\\
&=\frac{|\eta|^{2}}{2k}(e^{1}\wedge\bar{e}^{1}+\cdot\cdot\cdot e^{2k}\wedge\bar{e}^{2k}). \nonumber
 \end{align*}
It follows from (\ref{2.2}) that
\begin{align*}
	\Delta_{g}|\eta|_g^2=\sum\limits_{i=1}^n|D_i'\eta|^2+2k\sum\limits_{i=1}^nR_{i\bar{i}}\beta_{i\bar{i}}=\sum\limits_{i=1}^n|D_i'\eta|^2+|\eta|^2\sum\limits_{i=1}^{2k}R_{i\bar{i}}.
\end{align*}
and
 \begin{align*}
\langle\sqrt{-1}\partial\bar{\partial}\vert\eta\vert^2,\beta \rangle&=\sum\limits_{i=1}^{2k}\vert D'_i\eta\vert^2\cdot\beta_{i\bar{i}}+2k \sum\limits_{i=1}^{2k}R_{i\bar{i}j\bar{j}}\beta_{i\bar{i}}\beta_{j\bar{j}}\\
&=\frac{|\eta|^2}{2k}\sum\limits_{i=1}^{2k}|D_i'\eta|^2+\frac{|\eta|^4}{2k}\sum\limits_{i,j=1}^{2k}R_{i\bar{i}j\bar{j}}.
 \end{align*}
By formulas (\ref{2.4}) and (\ref{2.5}), we deduce that
 \begin{align*}
&\,\,\,\,\,\,\,\,\int_{|Z|=1,Z\in \sum_{2k}}\mathcal{C}_{a,b}(Z)d\theta(Z) \nonumber \\
&=a\int_{|Z|=1,Z\in \sum_{2k}}Ric(X,\overline{X})(Z)d\theta(Z)+b\int_{|Z|=1,Z\in \sum_{2k}}H(Z)d\theta(Z)) \nonumber \\
&=Vol(\mathbb{S}^{4k-1})[\frac{a}{2k}\sum_{i=1}^{2k}R_{i\overline{i}}+\frac{2b}{2k(2k+1)}\sum_{i,j=1}^{2k}R_{i\bar{i}j\bar{j}}]
 \end{align*}
where $\sum_{2k}$ is the $2k$-dimension subspace of $T_{x}^{1,0}M$ generated by $\{e_{i}\}_{i=1}^{2k}$. This implies that
\begin{align*}
\frac{a}{2k}\sum_{i=1}^{2k}R_{i\bar{i}}+\frac{2b}{2k(2k+1)}\sum_{i,j=1}^{2k}R_{i\bar{i}j\bar{j}}\geq C_{a,b}(x):=\min_{X\in T_x^{1,0}M\setminus\{0\}}\mathcal{C}_{a,b}(X).
 \end{align*}
Hence, if $a\geq0$ and $3a+2b\geq0$, we have
\begin{align*}
	&\frac{a}{2k}\Delta_{g}|\eta|_g^2\cdot|\eta|_g^2+\frac{b}{k(2k+1)}2k\langle\sqrt{-1}\partial\bar{\partial}\vert\eta\vert^2,\beta\rangle\\
	=&|\eta|_g^2\left(\frac{a}{2k}\sum\limits_{i=1}^n|D_i'\eta|^2+\frac{b}{k(2k+1)}\sum\limits_{i=1}^{2k}|D_i'\eta|^2\right)+|\eta|^4\left(\frac{a}{2k}\sum\limits_{i=1}^{2k}R_{i\bar{i}}+\frac{b}{k(2k+1)}\sum\limits_{i,j=1}^{2k}R_{i\bar{i}j\bar{j}}\right)\\
	\geq&C_{a,b}(x)\cdot|\eta|^4
\end{align*}
where it holds for any $x\in M$ such that $\eta\neq0$. Obviously, it also holds for any $x\in M$ such that $\eta=0$. From the formula $(\ref{2.1})$, we have
 \begin{align}
\int_{M}\Delta_{g}|\eta|_{g}^{2}\cdot|\eta|_{g}^{2}\cdot\omega^{n}&=2k\int_{M}\langle \sqrt{-1}\partial\bar{\partial}|\eta|^{2}_{g},\beta\rangle_{g}\cdot\omega^{n}\nonumber\\
&\,\,\,\,\,\,\,\,+(\sqrt{-1})^{(2k)^{2}}\frac{n!}{(n-2k-1)!}\int_{M}\sqrt{-1}\partial\bar{\partial}|\eta|_{g}^{2}\wedge\eta\wedge\bar{\eta}\wedge\omega^{n-2k-1}  \nonumber \\
&=2k\int_{M}\langle \sqrt{-1}\partial\bar{\partial}|\eta|^{2}_{g},\beta\rangle_{g}\cdot\omega^{n} \nonumber
 \end{align}
where the last equality holds from a well know result that any holomorphic $(p,0)$-form $\eta$ on a compact K\"{a}hler manifold is $d$-closed.
Now we have that
\begin{align*}
&(\frac{a}{2k}+\frac{b}{k(2k+1)})\int_{M}\Delta_{g}|\eta|_{g}^{2}\cdot|\eta|_{g}^{2}\cdot\omega^{n}\\
\geq &\frac{a}{2k}\int_{M}\Delta|\eta|^2\cdot|\eta|^{2}\cdot\omega^{n}
+\frac{b}{k(2k+1)}2k \int_{M}\langle\sqrt{-1}\partial\bar{\partial}|\eta|^2,\beta\rangle\omega^{n}\\
\geq &\int_M C_{a,b}(x)\cdot|\eta|^4\cdot\omega^n.
 \end{align*}
If $\mathcal{C}_{a,b}$ is nonnegative,  we have
   \begin{align*}
 -\int_{M}\left|d|\eta|^2\right|^2\cdot\omega^{n}=\int_{M}\Delta_{g}|\eta|_{g}^{2}\cdot|\eta|_{g}^{2}\cdot\omega^{n}\geq 0
 \end{align*}
Then $d|\eta|^{2}_{g}\equiv0$ on $M$, which implies that $|\eta|^{2}_{g}\equiv C$ for some constant $C>0$.
Since $C_{a,b}$ is
quasi-positive, we get that
 \begin{align*}
 0\geq C^2\int_{M}C_{a,b}(x)\cdot\omega^{n}>0
 \end{align*}
 which is a contradiction. Hence, there is no nonzero holomorphic $(2,0)$-form. The projectivity follows from Kodaira's theorem (\cite{Koda}, Theorem 1).
 \qed

\begin{rem}\label{rem1.3}
	If K\"{a}hler metric $g$ has quasi-positive $k$-Ricci curvature, then mixed curvature $\mathcal{C}_{k-1,n-k}$ is also quasi-positive by (\ref{2.3}). Hence, the Corollary \ref{1.3} holds. The proof for quasi-positive $2$-scalar curvature is essentially same. In fact, if K\"{a}hler metric $g$ has quasi-positive $2$-scalar curvature $S_{2}$, then $S_{k}=\sum\limits_{i,j=1}^kR_{i\bar{i}j\bar{j}}$ is also quasi-positive for
	$k\geq 2$, which is sufficient in the above arguement.
\end{rem}

\subsection{Rational connectedness}
In this section, we prove the the rational connectedness of Theorem 1.1. The proof shall follow closely Heier-Wong's method \cite{HW}, so we just briefly outline the proof.

It follows from (\ref{2.4}) and (\ref{2.5}) that
 \begin{align*}
S=\sum_{i,j=1}^{n}R_{i\overline{i}j\overline{j}}=\frac{n(n+1)}{(n+1)a+2b}\,\,\frac{1}{Vol(\mathbb{S}^{2n-1})}\int_{|Z|=1,Z\in T^{1,0}M}\mathcal{C}_{a,b}(Z)d\theta(Z)
 \end{align*}
where implies the scalar curvature $S$ is quasi-positive. By Heier-Wong's result \cite{HW1}, we know that $M$ is uniruled. In fact, we also have that
$$\int_Mc_1(K_M)\wedge\omega^{n-1}=	\frac{-1}{n\pi}\int_MS\omega^n<0.$$
 In (\cite{HW1}, Section 2), a linear algebra argument was given for the fact that any pseudo-effective line bundle $P$ on $M$ satisfies
$$\int_Mc_1(P)\wedge\omega^{n-1}\geq0.$$
From these inequalities, we easily deduce that $K_M$ must be not pseudo-effective and $M$ is uniruled by \cite{BDPP}. Therefore, $M$ admits a MRC fibration $f:M\rightarrow N$ from $M$ to a projective manifold $N$, and $K_N$ is pseudo-effective (\cite{Cam,KMM}).

We shall argue by contradiction. Assume that $m=\dim N\geq1$. Following the arguement as (Heier-Wong \cite{HW}, section 3.1 and section 3.2), there is a dense Zariski-open subset $Z'\subset N$ such that
$$f\vert_{f^{-1}(Z')}:f^{-1}(Z')\rightarrow Z'$$
is a holomorphic submersion, and the proof is reduced to show that for any fixed
small open set $V\subset Z'$, $$\int_{f^{-1}(V)}c_1(E)\wedge\omega^{n-1}=\int_{f^{-1}(V)}c_1(\Lambda^m E)\wedge\omega^{n-1}$$
is non-negative in general and positive for some $V$, where $E:=f^{*}T_N$. Here $E$ is the unique extension of $f^{*}T_N$ as a reflexive sheaf (see section 3.1 in \cite{HW}). In fact, this will contradict the fact that $K_N$ is pseudo-effective.

On $M\setminus f^{-1}(V)$, we have the standard exact sequence
$$0\rightarrow(T_{\left(M\setminus f^{-1}(V)\right)/\left(N\setminus V\right)},g_S)\rightarrow (T_M,g)\rightarrow (f^*T_N,g_E)\rightarrow0,$$
where $g_S$ and $g_E$ is the induced metric and quotient metrics respcetively. Then the curvature of $(E,g_E)$ is greater or equal the curvature of $(T_M,g)$.

For each $q\in V$, we write $M_{q}:=f^{-1}(q)$ for the fiber over $q$. For any fixed $p\in f^{-1}(V)$, we can choose an orthonormal basis $\{e_k\}_{k=1}^n$ of $T_p^{1,0}M$ such that $\{e_i\}_{i=1}^m$ is an orthonormal basis of $((T_pM_{f(p)})^\perp)_{p\in f^{-1}(V)}$. Note that
we use  the basis $\{e_i\}_{i=1}^{m}$ to denote the ``horizontal direction'' of the MRC fibration, while the ``horizontal direction'' is denoted by $\{e_{n-m+i}\}_{i=1}^{m}$ in section 3.2 of Heier-Wong \cite{HW}.
Then $(\Lambda^mE,G:=\Lambda^mg_E)$ satisfies that
$$\sqrt{-1}\Theta_{G}(u,\bar{u})\geq \sum\limits_{i=1}^mR_{u\bar{u}i\bar{i}},\ \forall u\in T_p^{1,0}M.$$
By a nontrivial calculation in (Heier-Wong, section 3.2) (also see (\cite{Mat2}, section 3)),
$$\int_{f^{-1}(V)}\sum\limits_{i=m+1}^n\sqrt{-1}\Theta_G(e_i,\bar{e}_i)\cdot\omega^n\geq0.
$$
By calculations similar to those in section 3.1 in \cite{HW} and Berger's lemma, for any small $V\subset Z'$, we can deduce that
\begin{align*}
&	\int_{f^{-1}(V)}c_1(\Lambda^m E)\wedge\omega^{n-1} \\
=&\int_{f^{-1}(V)}\sqrt{-1}\Theta_G\wedge\omega^{n-1}\\
	=&\int_{f^{-1}(V)}\sum\limits_{i=1}^n\sqrt{-1}\Theta_G(e_i,\bar{e}_i)\cdot\omega^n\\
	\geq&\frac{m(m+1)}{a(m+1)+2b}\int_{f^{-1}(V)}\left(\frac{a}{m}\sum\limits_{i=1}^n\Theta_G(e_i,\bar{e}_i)+\frac{2b}{m(m+1)}\sum\limits_{i=1}^m\Theta_G(e_i,\bar{e}_i)\right)\cdot\omega^n\\
	\geq&\frac{m(m+1)}{a(m+1)+2b}\int_{f^{-1}(V)}\left(\frac{a}{m}\sum\limits_{i=1}^mR_{i\bar{i}}+\frac{2b}{m(m+1)}\sum\limits_{i,j=1}^mR_{i\bar{i}j\bar{j}}\right)\cdot\omega^n\\
	\geq&\frac{m(m+1)}{a(m+1)+2b}\int_{f^{-1}(V)}C_{a,b}\cdot\omega^n\\
	\geq&0.
\end{align*}
Since $\{x\in M:C_{a,b}(x)>0\}$ is open, we can choose some $V_0$ such that $f^{-1}(V)\bigcap\{x\in M:C_{a,b}(x)>0\}\neq\emptyset$, then
$$\int_{f^{-1}(V_0)}c_1(\Lambda^m E)\wedge\omega^{n-1}\geq\frac{m(m+1)}{a(m+1)+2b}\int_{f^{-1}(V_0)}C_{a,b}\cdot\omega^n>0.$$
The proof is complete. \qed

\section{Real bisectional curvature}
For a compact Hermitian manifold $(M^{n},g)$, by Lemma 2.1 the following equalities hold:
\begin{align*}
-\int_{M}|\bar{\partial}|\eta|^{2}_{g}|^{2}_{g}\cdot\omega^{n}&=\int_{M}\Box_{g}|\eta|_{g}^{2}\cdot|\eta|_{g}^{2}\cdot\omega^{n} \nonumber \\
&=p\int_{M}\langle \sqrt{-1}\partial\bar{\partial}|\eta|^{2}_{g},\beta\rangle_{g}\cdot\omega^{n}\nonumber\\
&\,\,\,\,\,\,\,\,+(\sqrt{-1})^{p^{2}}\frac{n!}{(n-p-1)!}\int_{M}\sqrt{-1}\partial\bar{\partial}|\eta|_{g}^{2}\wedge\eta\wedge\bar{\eta}\wedge\omega^{n-p-1}
\end{align*}
where $\Box_{g}|\eta|_{g}^{2}$ is the complex Laplacian defined by $tr_{g}\sqrt{-1}\partial\bar{\partial}|\eta|_{g}^{2}$.
In non-K\"{a}hler case, it is difficult to determine the sign of the last term. But in some special case, we can estimate its sign.

\vspace{0.3cm}
\noindent {\bf Proof of Theorem 1.3:} We assume that $M$ admits a nonzero holomorphic $(2,0)$-form
$\eta$. If $g$ is Hermitian metric and complex dimension of the manifold is 3, then when $p=2$, we get that
\begin{align*}
	(\sqrt{-1})^{p^{2}}\frac{n!}{(n-p-1)!}&\int_{M}\sqrt{-1}\partial\bar{\partial}|\eta|_{g}^{2}\wedge\eta\wedge\bar{\eta}\wedge\omega^{n-p-1} \\
	&=3!\int_{M}\sqrt{-1}\partial\bar{\partial}|\eta|_{g}^{2}\wedge \eta\wedge\bar{\eta} \\
	&=3!\int_{M}\sqrt{-1}|\eta|_{g}^{2}\cdot(\partial\eta\wedge\overline{\partial\eta}) \\
	&\geq 0
\end{align*}
where these equalities follows from Stokes' theorem and $\bar{\partial}\eta=0$. If Hermitian metric $g$ has non-negative real bisectional curvature, then by (\ref{2.2}) we have
\begin{align*}
\langle\sqrt{-1}\partial\bar{\partial}\vert\eta\vert^2,\beta \rangle&=\sum\limits_{i=1}^{3}\vert D'_i\eta\vert^2\cdot\beta_{i\bar{i}}+2\sum\limits_{i=1}^{3}R_{i\bar{i}j\bar{j}}\beta_{i\bar{i}}\beta_{j\bar{j}}\geq 0
 \end{align*}
where $tr_{\omega}\beta=|\eta|_{g}^{2}$, $\beta=\sqrt{-1}\beta_{i\bar{i}}dz^{i}\wedge d\bar{z}^{i}$ at $x_{0}$ and $\beta_{i\bar{i}}\geq 0$.
Hence, we get that
\begin{align*}
0\geq-\int_{M}|\bar{\partial}|\eta|^{2}_{g}|^{2}_{g}\cdot\omega^{3}=\int_{M}\Box_{g}|\eta|_{g}^{2}\cdot|\eta|_{g}^{2}\cdot\omega^{3} \geq0
\end{align*}
Then $\bar{\partial}|\eta|^{2}_{g}\equiv0$ on $M$, which implies that $d|\eta|^{2}_{g}\equiv0$. So we have that $|\eta|^{2}_{g}\equiv C$ for some constant $C>0$. If real bisectional curvature is
quasi-positive, we get that
 \begin{align*}
 0\geq C^2\int_{M}\lambda(x)\cdot\omega^{3}>0
 \end{align*}
 where $\lambda(x)=\min B_{g}(e,a)$, which is a contradiction. Hence, there is no nonzero holomorphic $(2,0)$-form. The projectivity follows from Kodaira's theorem if
$M^{3}$ is K\"{a}hlerian.
 \qed

\vspace{0.3cm}
\noindent {\bf Proof of Theorem 1.4:} If $g$ is Hermitian metric, then when $p=2$, we get that
\begin{align*}
	&\,\,\,\,\,\,\,\,(\sqrt{-1})^{p^{2}}\frac{n!}{(n-p-1)!}\int_{M}\sqrt{-1}\partial\bar{\partial}|\eta|_{g}^{2}\wedge\eta\wedge\bar{\eta}\wedge\omega^{n-p-1} \\
	&=\frac{n!}{(n-3)!}\int_{M}\sqrt{-1}\partial\bar{\partial}|\eta|_{g}^{2}\wedge\eta\wedge\bar{\eta}\wedge\omega^{n-3} \\
	&=\int_{M}\sqrt{-1}|\eta|_{g}^{2}\partial \eta\wedge \overline{\partial \eta}\wedge\omega^{n-3}+\int_{M}\sqrt{-1}|\eta|_{g}^{2}\partial \eta\wedge \overline{ \eta}\wedge\bar{\partial}\omega^{n-3} \\
	&\,\,\,\,\,\,\,\,+\int_{M}\sqrt{-1}|\eta|_{g}^{2} \eta\wedge \overline{\partial \eta} \wedge \partial\omega^{n-3}
	+\int_{M}\sqrt{-1}|\eta|_{g}^{2} \eta\wedge \overline{\eta} \wedge \partial\bar{\partial}\omega^{n-3} \\
	&=0
\end{align*}
where since  $M$ is K\"{a}hlerian and $\partial\bar{\partial}\omega^{n-3}=0$. We can easily use the same arguments from Theorem \ref{1.6}  to complete the proof. \qed

\begin{rem}
Maybe it is interesting to consider that projectivity of conformally K\"{a}hler manifold with quasi-positive holomorphic sectional curvature.
\end{rem}

\end{document}